\documentclass[12pt]{amsart}
\usepackage{amsmath}
\usepackage{amssymb}
\usepackage{latexsym}

\newtheorem{theorem}{Theorem}[section]
\newtheorem{cor}[theorem]{Corollary}
\newtheorem{lemma}[theorem]{Lemma}

\newtheorem{prop}[theorem]{Proposition}
\newtheorem{defn}[theorem]{Definition}

\newenvironment{proof*}{\vskip 2mm\noindent {}}{\hfill $\Box$ \vskip 2mm}
\numberwithin{equation}{section}
\newcommand{\C}{{\mathbb{C}}}
\newcommand{\D}{{\mathbb{D}}}

\renewcommand{\P}{{\mathbb{P}}}

\newcommand{\Z}{{\mathbb{Z}}}
\newcommand{\I}{{\mathcal{I}}}

\renewcommand{\O}{{\mathcal{O}}}
\newcommand{\eps}{\varepsilon}

\begin{document}

\title[Limit of three-point Green functions]{Limit of three-point Green functions : the degenerate case
\footnote{2010 {\it Mathematics Subject Classification} 32U35, 32A27}
}

\author{DUONG Quang Hai, Pascal J. Thomas}

\address{Universit\'e de Toulouse\\ UPS, INSA, UT1, UTM \\
Institut de Math\'e\-matiques de Toulouse\\
F-31062 Toulouse, France} 
\email{quanghai@math.univ-toulouse.fr, pascal.thomas@math.univ-toulouse.fr}

\begin{abstract}
We investigate the limits of the ideals of holomorphic functions 
vanishing on three points in $\C^2$ when all three points tend to the
origin, and what happens to the associated pluricomplex Green functions.
This is a continuation of the work of Magnusson, Rashkovskii, Sigurdsson
and Thomas, where those questions were settled in a generic case. 
\end{abstract}

\keywords{pluricomplex Green function, complex Monge-Amp\`ere equation, ideals of holomorphic functions}

\thanks{}

\maketitle

\section{Introduction}
Let $\Omega$ be a hyperconvex bounded domain
in $\C^n$ containing the origin $0$ and let $\O(\Omega)$ denote the
space of holomorphic functions, respectively $PSH_-(\Omega)$ the space of
nonpositive plurisubharmonic functions on $\Omega$.  For every 
subset $S$ of $\Omega$ we let $\I(S)$ denote the ideal of all
holomorphic functions vanishing on $S$.  We  consider
ideals $\mathcal I$  such that their zero locus
$V(\mathcal I) := \{ z \in \Omega : f(z)=0, \forall f \in \mathcal I\}$ is a finite set.  Since the domain is pseudoconvex,
 there are finitely many global generators $\psi_j \in \mathcal O(\Omega)$ such that
for any $f\in \mathcal I$, there exists $h_j \in \mathcal O(\Omega)$ such that 
$f = \sum_j h_j \psi_j$, see e.g. \cite[Theorem 7.2.9, p. 190]{Ho}. 

\begin{defn}\cite{Ra-Si}
\label{greenideal}
Let $\mathcal I$ be an ideal of $\Omega$, and $\psi_j$ its generators. Then
$$
G_{\mathcal I}^\Omega (z) := \sup \big\{ u(z) : u \in PSH_-(\Omega), 
u(z) \le \max_j \log |\psi_{j}| + O(1) \big\}.
$$
\end{defn}

Note that the condition is meaningful only near $a \in V(\mathcal I)$.
In the special case when $S$ is a finite set in $\Omega$ and
$\mathcal I=\mathcal I(S)$, 
we write $G_{\mathcal I(S)}=G_S$. This case reduces to Pluricomplex Green functions with logarithmic singularities, already studied by many authors, 
e.g.~Demailly \cite{De1},  \cite{Za},  Lelong \cite{Lel}, and
Rashkovskii and Sigurdsson \cite{Ra-Si}.

Following the lead of \cite{MRST}, we want to study the limit of $G_{S_\eps}$
when $S_\eps$ is a set of points tending to the origin, and relate that to the
limit of the ideals $\I(S_\eps)$ (in a sense to be specified below, see \cite{MRST} for
more details). 

\begin{defn}
A (point based) ideal is a \emph{complete intersection ideal} if and only if 
it admits a set of $n$ generators, where $n$ is the dimension of the ambient space.
\end{defn}

The main result of \cite{MRST}, Theorem 1.11, states:

\begin{theorem}
\label{thmiff}
Let  $\mathcal I_\eps= \mathcal I(S_\eps)$, 
where $S_\eps$ is a set of $N$ points all tending to $0$ and
assume that
$\lim_{\eps\to 0} \mathcal I_\eps = \mathcal I$.
Then $(G_{\mathcal I_\eps})$ converges to $G_{\mathcal I}$ locally
uniformly on $\Omega \setminus \{0\}$ if and only if 
$\mathcal I$ is a complete intersection ideal.
\end{theorem}

Furthermore, \cite[Theorem 1.12, (i)]{MRST} works out the limits of 
Green functions when $N=3$ and the dimension $n=2$.  

Let $S_\eps:=\{a_1^\eps, a_2^\eps, a_3^\eps\}$.
For each pair of distinct indices, $i$, $j$, let $[a_i^\eps-a_j^\eps]=v_k^\eps \in \C\P^1$
where $\{ i,j,k\}=\{1,2,3\}$. The cases which are studied in \cite{MRST}
are those where there 
exist $i\neq j$ such that $\lim_{\eps\to 0} v_i^\eps$ and $\lim_{\eps\to 0} v_j^\eps$
exist and are distinct.  In those cases $\lim_{\eps\to 0} \I(S_\eps) = \frak M_0^2$
(the square of the maximal ideal at zero, i.e. the set of functions vanishing at zero
together with all their first derivatives), which is not a 
complete intersection ideal.  This sufficient condition is not necessary. 
We will give a characterization of the situations where the limit equals 
the square of the maximal ideal at zero.

The main goal of this note is to investigate the asymptotic behavior of ideals
and Green functions in 
the remaining (and most singular) case, when there exists $v \in \C^2$,
with $\|v\|=1$,
such that
\begin{equation}
\label{singhyp}
\lim_{\eps\to 0} v_i^\eps = [v] \mbox{ for }1\le i \le 3. 
\end{equation} 

We use the notation $z\cdot \bar w := z_1 \bar w_1 + z_2 \bar w_2$ for
$z, w \in \C^2$, and $\|z\|^2 := z \cdot \bar z$.

{\bf Numbering the points.} 

The notions we study do not depend on the order of the points in $S_\eps$, nor
does (\ref{singhyp}).
We choose an appropriate numbering. Set $d_i^\varepsilon := ||a_j^\varepsilon - a_k^\varepsilon||$, the Euclidean distances between two of the three points,
for $\{i, j, k\} = \{1, 2, 3\}, j, k \not= i$. For each $\eps$, number
the points so that $d_3^\varepsilon \geqslant d_1^\varepsilon \geqslant d_2^\varepsilon$. 

We perform a translation so that
$a_1^\varepsilon = (0,0)$.  Since the distance from $a_1^\varepsilon$ to the
origin tends to $0$ by hypothesis, this does not change any of the limits we
are studying, and we shall make this assumption henceforth.

Let $\theta$ be the acute angle between the lines
directed by $a_2^\varepsilon$ and $a_3^\varepsilon$, i.e. 
$\theta := \cos^{-1} \left( \frac{|a^\eps_2 \cdot \bar a^\eps_3 |}{\|a^\eps_2\| \|a^\eps_3\|}\right) $.

\begin{theorem}
\label{charmaxsq}
$\lim_{\eps\to 0} \I(S_\eps) = \frak M_0^2$ if and only if 
$\lim_{\eps\to 0} \frac{\|a^\eps_2\|}{|\theta|}=0$, or equivalently
\begin{equation}
\label{maxsq}
\lim_{\eps\to 0} \frac{\|a^\eps_2\|}{|\det (\frac{a^\eps_2}{\|a^\eps_2\|},\frac{a^\eps_3}{\|a^\eps_3\|})|}=0,
\end{equation}
where the determinant is taken with respect to an orthonormal basis.
\end{theorem}
The condition can be rephrased in a coordinate-free way using some elementary plane geometry: let $d$ be the diameter of the set $S_\eps$, let $\theta_1 \le \theta_2 \le \theta_3$ be the three angles (in $[0,\pi]$) determined by the triangle $a_1 a_2 a_3$; then 
we require that $d/\theta_2$ tend to $0$.

{\bf Choice of coordinates.} 

Now suppose that all three points converge to the origin along a common direction, 
that is, that 
there exists $v $ such that \eqref{singhyp} holds.
We reparametrize our family so that $|\eps|= \|a_1^\varepsilon - a_2^\varepsilon\|$.
We choose coordinates depending on $\eps$ 
so that $a_2^\varepsilon = (\varepsilon,0), a_3^\varepsilon = (\rho_1(\varepsilon), \rho_2(\varepsilon))$, where 
$\lim_{\eps\to 0} \rho_j(\varepsilon)=0$ for $j = 1, 2$, and 
$\lim_{\eps\to 0} \dfrac{\rho_2(\varepsilon)}{\rho_1(\varepsilon)} = 0$ (see
details before \eqref{newcoor}). 

Furthermore, $a_3^\varepsilon \in B(0; |\varepsilon|) \cap B(a_2^\varepsilon; |\varepsilon|)$ 
so we have 
$$|\rho_1(\varepsilon)| \leqslant \dfrac{1}{2}|\varepsilon|, |\rho_2(\varepsilon)| \leqslant \dfrac{\sqrt{3}}{2}|\varepsilon|. 
$$
Finally, we can write $a_3^\varepsilon = (\rho(\varepsilon), \delta(\varepsilon)\rho(\varepsilon))$, with $0 \not= |\rho(\varepsilon)| \leqslant \dfrac{1}{2}|\varepsilon|$, $\delta(\varepsilon) \rightarrow 0$.

\begin{theorem} 
\label{idealsconv}
Under the above hypotheses:

(1) \; If $\underset{\varepsilon \longrightarrow 0}{\lim} \dfrac{\delta(\varepsilon)}{\rho(\varepsilon) - \varepsilon} = m \not=  \infty$, then 
$$
\underset{\varepsilon \longrightarrow 0}{\lim} \mathcal I_\varepsilon = \mathcal I = \left\langle z_2 - m z_1^2, z_1^3\right\rangle,
$$
so by Theorem \ref{thmiff},
$$
\underset{\varepsilon \longrightarrow 0}{\lim} G_\varepsilon(z) = G_{\mathcal I}(z) = \max\big(\log|z_2 - m z_1^2|, 3\log|z_1|\big).
$$

(2) \; If $\underset{\varepsilon \longrightarrow 0}{\lim} \dfrac{\delta(\varepsilon)}{\rho(\varepsilon) - \varepsilon} = \infty$ (equivalently $\underset{\varepsilon \longrightarrow 0}{\lim} \dfrac{\delta(\varepsilon)}{\varepsilon} = \infty$),
then 
$$\underset{\varepsilon \longrightarrow 0}{\lim} \mathcal I_\varepsilon = \mathfrak{M}_0^2,
$$
so by Theorem \ref{thmiff}, the Green function cannot converge to $G_{ \mathfrak{M}_0^2}(z)=2 \log \|z\| + O(1)$. 
\end{theorem} 

We still need to understand to which limit the Green function may converge in case (2), at 
least in the model case where $\Omega = \D^2$. 
Unfortunately, we could only get some partial estimates.

\begin{prop}
\label{Greenest} 
(1) For any $z = (z_1, z_2) \in \mathbb{D}^2 \backslash \{(0,0)\}$, 
\begin{equation*} 
\underset{\varepsilon \longrightarrow 0}{\lim\inf} \; G_\varepsilon(z_1, z_2) \geqslant \max\big(2 \log|z_1|,\dfrac{3}{2}\log |z_2|\big).
\end{equation*}

(2) For any $z = (z_1, z_2) \in \mathbb{D}^2 \backslash \{(0,0)\}$ such that $|z_2| \leqslant |z_1|^2$, 
$$
\underset{\varepsilon \longrightarrow 0}{\lim} \; G_\varepsilon(z_1, z_2) = 2\log |z_1|.
$$
\end{prop}

In general, it is more difficult to get upper than lower bounds on the limits of Green functions.  We only could get a general upper bound under rather special hypotheses on 
the configuration, roughly speaking that the angle with vertex at the origin formed by the three points should tend to $0$ very slowly.

\begin{prop}
\label{deltabig} 
If $\dfrac{\log |\delta|}{\log |\varepsilon|} \longrightarrow 0$ as $\varepsilon \rightarrow 0$, then 
for any $z = (z_1, z_2) \in \mathbb{D}^2 \backslash \{(0,0)\}$
\begin{equation*} 
\underset{\varepsilon \longrightarrow 0}{\lim\sup} \; G_\varepsilon(z) \leqslant 
\dfrac{3}{2} \log\max(|z_1|, |z_2|).
\end{equation*}
\end{prop}

\begin{cor}
\label{partialconv} 
Under the hypotheses of the previous Proposition, 
if $|z_2|\ge |z_1|$, then $\lim_{\eps\to0} G_\varepsilon(z) = \dfrac{3}{2} \log |z_2|$.
\end{cor}


\section{Proof of Theorems \ref{charmaxsq} and \ref{idealsconv}}

\subsection{Preliminary facts.}

First we need a notion of convergence of ideals, inspired by Hausdorff convergence. 
This is taken from \cite{MRST}.
Let $\Omega$ be a bounded pseudoconvex domain in $\C^n$. Let $E \subset \C$ such that $\Bar{E} \ni 0$ the set of parameters along which we take limits.  Convergence of holomorphic functions is always understood uniformly on compacta.

\begin{defn}  If $(\mathcal{I}_\varepsilon)_{\varepsilon \in E}$ are ideals
in $\O(\Omega)$, we define  
\begin{equation*} 
\begin{aligned}
\liminf\limits_{E\ni\eps \to 0}\I_\eps & := 
\{f \in \O(\Omega) : \; \mbox{for all}\; \varepsilon \in E, \exists f_\varepsilon \in \mathcal{I}_\varepsilon, f_\varepsilon \longrightarrow f\\
&\mbox{when}\; \varepsilon \longrightarrow 0\}. 
\end{aligned}
\end{equation*}
Likewise
$\limsup\limits_{E\ni\eps \to 0}\I_\eps $ is the vector space generated by
\begin{equation*} 
\begin{aligned}
\{f \in &\O(\Omega) :  \exists (\varepsilon_j)_{j \in \Z_+} \subset E, \varepsilon_j \longrightarrow 0 \; \mbox{when}\; j \longrightarrow \infty; \exists f_j \in \mathcal{I}_{\varepsilon_j},\\
&f_j \longrightarrow f\; \mbox{when}\; j \longrightarrow \infty\}.
\end{aligned}
\end{equation*}
We say that $(\mathcal{I}_\varepsilon)_{\varepsilon \in E}$ converges to $\I$ 
if and only if $\liminf\limits_{E\ni\eps \to 0}\I_\eps  = \limsup\limits_{E\ni\eps \to 0}\I_\eps  = \I$, and write $\underset{\varepsilon \longrightarrow 0, \varepsilon \in E}{\lim} \mathcal{I}_\varepsilon = \I$.
\end{defn}

Of course $\liminf\limits_{E\ni\eps \to 0}\I_\eps 
\subset 
\limsup\limits_{E\ni\eps \to 0}\I_\eps $, and they are both ideals.

Denote the Taylor expansion and Taylor polynomial of a holomorphic function $f$ by
\begin{equation*} 
f(z) = f(z_1, z_2) = \sum_{j,k=0}^\infty a_{jk} z_1^j z_2^k ; \quad
P_m(f)(z) := \underset{j+k \leqslant m}{\sum_{j,k}} a_{jk} z_1^j z_2^k.
\end{equation*}

It follows from the Cauchy formula on the distinguished boundary of $\mathbb{D}^2$ that
\begin{lemma} 
\label{firstterms}
 Let $m \in \mathbb{N}^\ast$, $U$ a bidisk centered at $(0,0)$, relatively compact in $\mathbb{D}^2$. There exists $C = C(m, U)$ such that for any $f\in \O(\mathbb{D}^2)$ with $\underset{\mathbb{D}^2}{\sup}||f||  \leqslant 1$, there exist holomorphic functions $r_{j,k} \in \O(\mathbb{D}^2)$ satisfying : for $j + k = m+1$, $\underset{U}{\sup}|r_{j,k}| \leqslant C, 0 \leqslant j \leqslant m + 1$ and for $z = (z_1, z_2) \in U$, then
\begin{equation*} 
f(z) = P_m(f)(z) + R_{m+1}(z) = P_m(f)(z) +  \sum_{j=0}^{m+1} r_{j,m+1-j}(z) z_1^j z_2^{m+1-j}.
\end{equation*}
\end{lemma}

\subsection{Proof of the sufficiency in Theorem \ref{charmaxsq}.}

Suppose that $f \in \limsup_\eps \mathcal I_\eps$. 
This means that there exists some subset $E \subset \C$
such that $0 \in \overline E$ and 
 a family of holomorphic functions  $\{f^\varepsilon, \eps \in E\}$,
with $f^\varepsilon \in \mathcal{I}_\varepsilon$, $\eps \in E$,
converging to $f$ uniformly 
on a fixed neighborhood $U$ of the origin.
Observe that all the Taylor coefficients will have to converge.  

Since $f^\varepsilon(a_1^\varepsilon) = 0$, $a_{0,0}^\varepsilon = 0$ for any $\varepsilon$. 
Applying Lemma \ref{firstterms} for $m = 1$, if $U \Subset U' \Subset \Omega$,
\begin{equation*} 
f^\varepsilon(z_1, z_2) = a_{1,0}^\varepsilon z_1 + a_{0,1}^\varepsilon z_2 + R_2(z_1, z_2)
\end{equation*}
with $|R_2(z_1, z_2)| \le C \|z\|^2$, where $C$ only depends on
$U$, $U'$ and $\sup_{U'}|f^\varepsilon|$.

Applying this to $z=a_i^\eps$, dividing by $\|a^\eps_i\|$
and writing $\nabla f^\eps (0) := (a_{1,0}^\varepsilon,a_{0,1}^\varepsilon)$,
we find
$$
\frac{a^\eps_i}{\|a^\eps_i\|} \cdot \nabla f^\eps (0) = O ( \|a^\eps_i\| ), \quad i=2,3.
$$
Write $M$ for the $2 \times 2$ matrix with rows given by 
the coordinates of $\frac{a^\eps_2}{\|a^\eps_2\|}$
and $\frac{a^\eps_3}{\|a^\eps_3\|}$.  Then $\|M\|=O(1)$ and 
$\|M^{-1}\|=O\left(|\det (\frac{a^\eps_2}{\|a^\eps_2\|},\frac{a^\eps_3}{\|a^\eps_3\|})|^{-1}\right)$.
Since by our choice of numbering, $\|a^\eps_3\|\le \|a^\eps_2\|$, we have
$$
\nabla f^\eps (0) = O ( \|M^{-1}\| \|a_2^\varepsilon\| ),
$$
so that if condition \eqref{maxsq} is met, then 
$\lim_{E\ni \eps \to 0} \nabla f^\eps (0) =0$, thus $f \in \frak M_0^2$.  
We have proved that condition  \eqref{maxsq} implies
that $ \limsup_\eps \mathcal I_\eps \subset \frak M_0^2$. 

To prove the inclusion $\frak M_0^2 \subset \liminf_\eps \mathcal I_\eps $,
it will be easier to take suitable coordinates.

For each pair of distinct indices, $i$, $j$, let $a_i^\eps-a_j^\eps=u_k^\varepsilon \in \mathbb{C}^2$ where $\{i, j, k\} = \{1, 2, 3\}$. For $|\eps|$ small enough,
 $\theta \not= 0$ so $\{u_3^\varepsilon, u_2^\varepsilon\}$ are linearly independent. 
Using the Gram-Schmidt orthogonalization process, 
we get $\mathfrak{B}_\varepsilon := \{e_1^\varepsilon, e_1^\varepsilon\}$ 
an orthonormal basis of $\mathbb{C}^2$, where 
$e_j^\varepsilon = \dfrac{\upsilon_j^\varepsilon}{||\upsilon_j^\varepsilon||}$, for $j= 1, 2$, 
$\upsilon_1^\varepsilon := u_3^\varepsilon$ and $\upsilon_2^\varepsilon := u_2^\varepsilon - \dfrac{u_2^\varepsilon.\overline{u_3^\varepsilon}}{||u_3^\varepsilon||^2}.u_3^\varepsilon$. 
If $z = (z_1, z_2) \in \Omega$ its new coordinates $(z_1^\varepsilon, z_2^\varepsilon)$ 
with respect to $\mathfrak{B}_\varepsilon$ are given by
\begin{equation}
\label{newcoor}
z_1^\varepsilon = z \cdot \bar e_1^\varepsilon , 
z_2^\varepsilon = z \cdot \bar e_2^\varepsilon .
\end{equation}

Denote the coordinates of the points in this new basis as before.
 Theorem \ref{idealsconv}
(condition \eqref{maxsq} implies that
the angle between $v_2^\eps$ and $v_3^\eps$ tends to $0$,
but no convergence of the basis is needed).  Then $|\delta(\eps)|=\tan \theta$,
so our hypothesis implies that $\lim \eps/\delta= \lim (\rho(\eps)-\eps)/\delta=0$. 

The following polynomials are in $\mathcal{I}_\varepsilon$:
\begin{equation*}
\begin{aligned}
Q_1^\varepsilon(z) &= (z_1^\varepsilon)^2 - \eps z_1^\varepsilon - \dfrac{\rho - \varepsilon}{\delta}z_2^\varepsilon;\\
Q_2^\varepsilon(z) &= z_2^\varepsilon \big(z_1^\varepsilon - \rho\big);\\
Q_3^\varepsilon(z) &= z_2^\varepsilon \big(z_2^\varepsilon - \delta \rho\big).
\end{aligned}
\end{equation*}
Let $\alpha_{ij} := e_j^\varepsilon \cdot \bar e_i $, for $1 \leqslant i, j \leqslant 2$,
so that
\begin{eqnarray*}
z_1 &= \alpha_{11} z_1^\varepsilon + \alpha_{12} z_2^\varepsilon , \\
z_2 &= \alpha_{21} z_1^\varepsilon + \alpha_{22} z_2^\varepsilon .
\end{eqnarray*}
If we let 
\begin{equation*}
\begin{aligned}
f_1^\varepsilon (z) &= \alpha_{11}^2 Q_1^\varepsilon(z) + 2\alpha_{11}\alpha_{12}
Q_2^\varepsilon(z) + \alpha_{12}^2Q_3^\varepsilon(z),\\
f_2^\varepsilon (z) &= \big(\alpha_{11}\alpha_{22}+\alpha_{12}\alpha_{21}\big)Q_2^\varepsilon(z) +  \alpha_{11} \alpha_{21}Q_1^\varepsilon(z) + \alpha_{12} \alpha_{22}Q_2^\varepsilon(z), \\
f_3^\varepsilon (z) &= \alpha_{21}^2Q_1^\varepsilon(z) + 2\alpha_{21}\alpha_{22}Q_2^\varepsilon(z) + \alpha_{22}^2Q_3^\varepsilon(z)  ,\\
\end{aligned}
\end{equation*}
then 
\begin{equation*}
\begin{aligned}
z_1^2 &= \underset{\varepsilon \longrightarrow\; 0}{\lim}\; f_1^\varepsilon (z) \in \underset{\varepsilon \longrightarrow\; 0}{\lim\inf}\;\mathcal{I}_\varepsilon;\\
z_1z_2 &= \underset{\varepsilon \longrightarrow\; 0}{\lim}\; f_2^\varepsilon (z) \in \underset{\varepsilon \longrightarrow\; 0}{\lim\inf}\;\mathcal{I}_\varepsilon;\\
z_2^2 &= \underset{\varepsilon \longrightarrow\; 0}{\lim}\; f_3^\varepsilon (z) \in \underset{\varepsilon \longrightarrow\; 0}{\lim\inf}\;\mathcal{I}_\varepsilon,
\end{aligned}
\end{equation*}
which proves that $\frak M_0^2 \subset \liminf_\eps \mathcal I_\eps $, and thus that
condition \eqref{maxsq} is sufficient for the claimed convergence. 

\subsection{Proof of Theorem \ref{idealsconv}.}

Statement (2) in Theorem \ref{idealsconv} is a special case of the proof above.
We now turn to statement (1). 

We slightly modify the varying basis $\mathfrak{B}_\varepsilon$ from the previous proof.
The hypothesis of Theorem \ref{idealsconv} implies that  
$\lim_{\eps} [e_1^\eps]$ exists, so multiplying $e_1^\eps, e_2^\eps$ by appropriate complex
numbers of modulus one, we get a basis 
$\tilde {\mathfrak{B}}_\varepsilon = (\tilde e_1^\eps ,\tilde e_2^\eps )$ such that 
$\lim_{\eps} \tilde e_1^\eps =:e_1$ and $\lim_{\eps} \tilde e_2^\eps =:e_2$ exist.  
We denote by
$(z_1^\eps, z_2^\eps)$, respectively $(z_1, z_2)$, the coordinates in $(\tilde e_1^\eps ,\tilde e_2^\eps )$, resp. $(e_1,e_2)$. 

Finally, given any function $f$ expressed in the $(z_1, z_2)$-coordinates, 
we denote by $\tilde f$ the function computed in the $(z_1^\eps, z_2^\eps)$-coordinates
(i.e. if $z$ and $\tilde z$ are the coordinates of the same point, $f(z)=\tilde f(\tilde z)$).  We write
$$
\tilde f^\eps (z_1^\eps, z_2^\eps)= \sum_{j,k} \tilde a^\eps_{ij} (z_1^\eps)^j  (z_2^\eps)^k
$$
(both the function and the coordinates depend on $\eps$). 
Since this is a linear change of variables, the various Taylor coefficients are 
obtained from the chain rule by linear formulae, and since the change of variables
matrix tends to the identity as $\eps \to 0$, $\lim_\eps a^\eps_{ij} = \lim_\eps \tilde a^\eps_{ij}$ if the latter exists.

Again, let $f \in \liminf_\eps \mathcal I_\eps$, i.e.
$f=\lim_\eps f^\varepsilon$ with uniform convergence 
on a fixed neighborhood of the origin, $f^\varepsilon \in \mathcal I_\eps$.  
 Applying  Lemma \ref{firstterms} for $m = 2$, taking $a_{0,0}^\varepsilon=0$
 into account,
\begin{equation*} 
\tilde f^\varepsilon(z_1, z_2) = \tilde a_{1,0}^\varepsilon z_1 + \tilde a_{0,1}^\varepsilon z_2 + \tilde a_{2,0}^\varepsilon z_1^2 + \tilde a_{0,2}^\varepsilon z_2^2 + \tilde a_{1,1}^\varepsilon z_1z_2 + R_3(z_1, z_2)
\end{equation*}
with $|R_3(z_1, z_2)| \le C \|z\|^3$.

Since $\tilde f^\varepsilon(a_2^\varepsilon) = \tilde f^\varepsilon(\varepsilon,0) = 0$ we have $\tilde a_{1,0}^\varepsilon \varepsilon + \tilde a_{2,0}^\varepsilon \varepsilon^2 + R_3(\varepsilon,0) = 0$. Thus
\begin{equation}
\label{equa1}
\tilde a_{1,0}^\varepsilon = - \tilde a_{2,0}^\varepsilon \varepsilon - \dfrac{R_3(\varepsilon,0)}{\varepsilon}
\end{equation}
for any $\varepsilon$. 

Thus $\dfrac{\partial \tilde f}{\partial z_1}(0,0) = \tilde a_{1,0} =  
\lim_{\eps\to0} \tilde a_{1,0}^\varepsilon = 0$. 

Furthermore, setting $\rho = \rho(\varepsilon), \delta = \delta(\varepsilon)$, 
from
$\tilde f^\varepsilon(a_3^\varepsilon) = \tilde f^\varepsilon(\rho, \delta\rho) = 0$, and  \eqref{equa1}
we deduce
\begin{equation} 
\label{equa2}
\begin{aligned}
&\bigg[- \tilde a_{2,0}^\varepsilon \varepsilon - \dfrac{R_3(\varepsilon,0)}{\varepsilon}\bigg]\rho + \tilde a_{0,1}^\varepsilon\delta\rho + \tilde a_{2,0}^\varepsilon \rho^2 + 
\tilde a_{0,2}^\varepsilon \delta^2\rho^2 \\
&+ \tilde a_{1,1}^\varepsilon \delta\rho^2 + R_3(\rho, \delta\rho) = 0,
\end{aligned}
\end{equation}
and dividing by $\rho(\rho-\varepsilon)$,
\begin{equation*} 
\label{equa3}
\begin{aligned}
\tilde a_{2,0}^\varepsilon + \tilde a_{0,1}^\varepsilon \dfrac{\delta}{\rho-\varepsilon} + \tilde a_{0,2}^\varepsilon \delta^2 \dfrac{\rho}{\rho-\varepsilon} + \tilde a_{1,1}^\varepsilon \dfrac{\delta}{\rho-\varepsilon}\rho + \dfrac{R_3(\rho, \delta\rho)}{\rho(\rho-\varepsilon)} 
- \dfrac{R_3(\varepsilon,0)}{\eps(\rho-\varepsilon)}= 0.
\end{aligned}
\end{equation*}
Note that since $|\rho(\varepsilon)| \leqslant \dfrac{1}{2}|\varepsilon|$, 
 $\dfrac{2}{3} \leqslant \bigg| \dfrac{\varepsilon}{\rho-\varepsilon}\bigg| \leqslant 2$, 
and $\bigg|\dfrac{\rho}{\rho-\varepsilon}\bigg| = \bigg|\dfrac{\varepsilon}{\rho-\varepsilon}\bigg| \bigg|\dfrac{\rho}{\varepsilon}\bigg| \leqslant 2\bigg|\dfrac{\rho}{\varepsilon}\bigg| \leqslant 1$. 

Since $R_3(\rho, \delta\rho) = O(\rho^3)$, $\underset{\varepsilon \longrightarrow 0}{\lim} \dfrac{R_3(\rho, \delta\rho)}{\rho(\rho-\varepsilon)} = 0$ and passing to
the limit as explained above,
$a_{2,0} + m a_{0,1} = 0$.  Thus
\begin{equation*} 
\begin{aligned}
{}&\underset{\varepsilon \longrightarrow 0}{\lim\sup}\; \I_\varepsilon \subset \I := \bigg\{f \in \mathcal{O}(\mathbb{D}^2) : f(0,0) =  \dfrac{\partial f}{\partial z_1}(0,0) = 0   \mbox{ and }\\
 &\dfrac{1}{2}\dfrac{\partial^2 f}{\partial z_1^2}(0,0) + m \dfrac{\partial f}{\partial z_2}(0,0) = 0\bigg\}
= \left\langle z_2 - m z_1^2, z_2^2, z_1z_2, z_1^3\right\rangle.
\end{aligned}
\end{equation*}
Since  $z_1z_2 = z_1(z_2 - m z_1^2) + m z_1^3$ et $z_2^2 = (z_2+ m z_1^2)(z_2 - m z_1^2) + (m^2 z_1)z_1^3$, we have $z_1z_2, z_2^2 \in \left\langle z_2 - m z_1^2, z_1^3\right\rangle$. 
Thus $\I = \left\langle z_2 - m z_1^2, z_1^3\right\rangle$. 

Conversely,
\begin{equation*} 
z_2 - m z_1^2 = \underset{\varepsilon \longrightarrow 0}{\lim} \bigg(z^\eps_2 - \dfrac{\delta}{\rho-\varepsilon}z^\eps_1(z^\eps_1-\varepsilon)\bigg) \in \underset{\varepsilon \longrightarrow 0}{\lim\inf}\; \I_\varepsilon, \mbox{ and}
\end{equation*}
\begin{equation*} 
z_1^3 = \underset{\varepsilon \longrightarrow 0}{\lim} z^\eps_1(z^\eps_1-\varepsilon)\big(z^\eps_1-\rho(\varepsilon)\big) \in \underset{\varepsilon \longrightarrow 0}{\lim\inf}\; \I_\varepsilon.
\end{equation*}
Thus $\I \subset \underset{\varepsilon \longrightarrow 0}{\lim\inf}\; \I_\varepsilon$. 
We have proved $\underset{\varepsilon \longrightarrow 0}{\lim}\; \I_\varepsilon = \I = \left\langle z_2 - mz_1^2, z_1^3\right\rangle$.

Since $\I$ admits a representation by two generators, the second statement in (1) 
follows from Theorem \ref{thmiff}.

\subsection{Proof of the necessity in Theorem \ref{charmaxsq}.}

If $\frac{\|a^\eps_2\|}{|\theta|}$ does not tend to $0$, 
we can find a sequence
$\eps_j \to 0$ such that $\frac{\|a^{\eps_j}_2\|}{|\theta|} \geqslant c >0$
and therefore along this sequence $\theta$ tends to $0$, so the distances in $\C\P^1$
between all $[a^\eps_i-a^\eps_j]$ tend to $0$.  Also, $\theta \sim \tan \theta
\sim \sin \theta$. 
Passing to a further subsequence, we may assume that
$a^\eps_2/\|a^\eps_2\|$ converges. Using
the coordinates and notations of Theorem \ref{idealsconv},
we have that $\dfrac{\delta(\varepsilon_j)}{\rho(\varepsilon_j) - \varepsilon_j}$
is a bounded sequence in $\C$.  
Then, passing to another subsequence, we may assume
that ${\lim_j} \dfrac{\delta(\varepsilon_j)}{\rho(\varepsilon_j) - \varepsilon_j} = m \in \C$
 and so we are in the situation of Theorem \ref{idealsconv}, statement (1). So the limit of the ideals $\mathcal I_\eps$
 along this subsequence 
contains the function $z_2$ (given by this appropriate coordinate system),
which implies that $\limsup_\eps \mathcal I_\eps \not\subset \frak M_0^2$. 

We remark that in this case, if $[a^\eps_2]$ does not converge, then we can find two 
different limit values for it, and two different functions of degree $1$ in 
$\limsup_\eps \mathcal I_\eps$, so that $\limsup_\eps \mathcal I_\eps =\frak M_0$
and (for reasons of length) the family $( \mathcal I_\eps)$ cannot converge
to any limit ideal.  This is in contrast to the other case, where no convergence
of the varying basis $\frak B_\eps$ was required. 

\section{Proof of Proposition  \ref{Greenest} }

The definition of the Green function implies that for any $f\in \O(\D^2)$, with
$\|f\|_\infty \le 1$, such that $f(a_j^\eps)=0$, $1\le j \le 3$, then
$G_\eps (z) \ge \log|f(z)|$.

First consider the polynomial
$Q(z_1,z_2) = -\varepsilon z_1 + z_1^2 + \dfrac{\varepsilon-\rho}{\delta}z_2$. Define
\begin{equation*} 
Q_1(z_1, z_2) := \bigg(\underset{\underset{|w_2| < 1}{|w_1| < 1}}{\sup}\big|Q(w_1,w_2)\big|\bigg)^{-1} Q(z_1,z_2).
\end{equation*}
It is easy to see that it satisfies the conditions above, so letting
 $\varepsilon$ tend to $0$,
\begin{equation}\label{equa4}
\underset{\varepsilon \longrightarrow 0}{\lim\inf} \; G_\varepsilon(z_1, z_2) \geqslant \log |z_1|^2 = 2\log |z_1|,
\end{equation}
for any $z = (z_1, z_2) \in \mathbb{D}^2 \backslash \{(0,0)\}$. 

To get the other part of the estimate, consider the three lines passing through two points $a_i^\varepsilon$ and $a_j^\varepsilon$, $(i \not= j, i, j \in \{1, 2, 3\})$, with 
the following equations :
$$
l_1^\varepsilon(z) = z_2 ;\quad
 l_2^\varepsilon(z) = z_2 - \delta(\varepsilon)z_1;\quad
  l_3^\varepsilon(z) = z_2 - \delta(\varepsilon)\dfrac{\rho(\varepsilon)}{\rho(\varepsilon) - \varepsilon}(z_1 - \varepsilon). 
$$
Since any pole belong to two of the lines,
 $P(z_1,z_2) := l_1^\varepsilon(z) \cdot l_2^\varepsilon(z) \cdot l_3^\varepsilon(z)
 \in \I_\varepsilon^2$. So
 \begin{equation*} 
G_\varepsilon(z_1, z_2) \geqslant \dfrac{1}{2} \log \dfrac{|P(z_1,z_2)|}{||P||_\infty}.
\end{equation*}
 Furthermore
 $$
 ||P||_\infty \leqslant 1(1+|\delta|)(1+|\delta|.(1+|\varepsilon|)), 
 $$
so $\lim_{\eps\to0}||P||_\infty =1$. Letting $\varepsilon$ tend to $0$, 
\begin{equation*} 
\underset{\varepsilon \longrightarrow 0}{\lim\inf} \; G_\varepsilon(z_1, z_2) \geqslant \dfrac{3}{2}\log |z_2|,
\end{equation*}
for any $z = (z_1, z_2) \in \mathbb{D}^2 \backslash \{(0,0)\}$. 

To obtain Part (2) of the Proposition, we observe as in \cite{MRST} that 
$$
G_\eps (z) \le G_{\{a_1^\eps, a_2^\eps\}}(z) = 
\max \left( \log \left|z_1 \frac{\eps-z_1}{1-\bar \eps z_1} \right| , \log |z_2|\right),
$$
and this last function tends to $2 \log |z_1|$ when $|z_2|\le |z_1|^2$.
 \hfill $\square$\\

\section{Proof of Proposition \ref{deltabig}}

We will estimate the function by restricting it to well-chosen families of analytic
disks.  By Proposition \ref{Greenest} (2), there is no loss in assuming $z_2\neq 0$. 

Let $Z_1 := \dfrac{z_1-\varepsilon}{\|z\|_\infty}r_\varepsilon$, $Z_2 := \dfrac{z_2}{\|z\|_\infty}r_\varepsilon$, where $r_\eps:= 1- \dfrac{|\eps||z_1|}{|z_1|-|\varepsilon|} \le 1$. 
If $z_1=1$, $r_\eps=1$; if $z_1 \neq 0$,  $r_\eps$ is defined for
  $|\varepsilon| < |z_1|$ and
in this case $0 < r_\varepsilon \rightarrow 1$ as $\varepsilon \rightarrow 0$.  
Note that $|Z_1|, |Z_2| \le 1$.
Set
 \begin{equation*} 
\Psi(\zeta) := (\varepsilon+Z_1\zeta, Z_2\zeta) = \bigg(\varepsilon + \dfrac{z_1-\varepsilon}{\|z\|_\infty}r_\varepsilon\zeta, \dfrac{z_2}{\|z\|_\infty}r_\varepsilon\zeta\bigg).
\end{equation*}

For $\zeta \in \mathbb{D}$, we have 
$\bigg|\dfrac{z_2}{\|z\|_\infty} r_\varepsilon\zeta\bigg| < r_\varepsilon \le 1$ and $\bigg|\varepsilon+\dfrac{z_1-\varepsilon}{\|z\|_\infty}r_\varepsilon\zeta\bigg| < 1$. Since, when $z_1\neq 0$, 
\begin{equation*} 
\bigg| r_\varepsilon\zeta + \dfrac{\varepsilon\|z\|_\infty}{z_1 - \varepsilon}\bigg| 
< |\zeta|\left(  1- \dfrac{|\eps||z_1|}{|z_1|-|\varepsilon|} \right) + 
\dfrac{|\varepsilon| |z_1|}{|z_1| - |\varepsilon|} <1.
\end{equation*}
Thus $\Psi(\zeta) \in \mathbb{D}^2$. Furthermore, for $|\varepsilon| < \varepsilon_0$, 
we have $\bigg|\dfrac{\|z\|_\infty}{r_\varepsilon} \bigg| < 1$ et $\Psi\big(\dfrac{\|z\|_\infty}{r_\varepsilon}\big) = (z_1,z_2)$, $\Psi(0) = a_1^\varepsilon = (\varepsilon,0)$.

Let $u \in PSH_{\_}(\mathbb{D}^2)$ a function in the defining family of $G_\varepsilon$.
Set $u_2 := u \circ \Psi$. For $\zeta \in \mathbb{D}$, 
$$
u_2(\zeta) =  u\big(\varepsilon+Z_1\zeta, Z_2\zeta\big) \leqslant \log\max \big(|Z_1\zeta|, |Z_2\zeta|\big)+O(1) \leqslant \log |\zeta|+O(1).
$$
So write $u_3(\zeta) := u_2(\zeta) - \log |\zeta|$, for $\zeta \in \mathbb{D} \backslash \{0\}$. Since $u_2(\zeta) \in SH_{\_}\big(\mathbb{D} \backslash \{0\}\big)$ 
and $\log |\zeta| \in H\big(\mathbb{D} \backslash \{0\}\big)$,  $u_3 \in SH_{\_}\big(\mathbb{D} \backslash \{0\}\big)$. Near $0$, $u_3(\zeta) = u_2(\zeta) - \log |\zeta| \leqslant \log |\zeta|+O(1) - \log |\zeta| = O(1)$, so $u_3$ is bounded
in a neighborhood of $0$. By the removable singularity theorem for 
subharmonic functions  \cite[Theorem 3.6.1]{Rans}, we can extend
$u_3$ to a function in $SH_{\_}\big(\mathbb{D}\big)$.

\begin{lemma}
\label{smalldisk}
For any $u\in PSH_{\_}(\mathbb{D}^2)$ a function in the defining family of $G_\varepsilon$, 
there exist $z_0 \in \D$ and constants $C_1, C_4>0$ such that 
for any $\zeta \in D_0:=D(z_0, C_1 |\eps|^2)$, 
$$
u_3 (\zeta) \le \log \left| \frac\eps\delta \right| + C_4.
$$
\end{lemma}

Let us postpone the proof of this lemma (which will use the other two poles and another family of analytic discs).  Notice that it doesn't use the hypothesis 
$\underset{\varepsilon \longrightarrow 0}{\lim} \dfrac{\log |\delta|}{\log|\varepsilon|} = 0$.

For any $\xi \in \mathbb{D}$, set $u_4(\xi) := u_3 \circ \varPhi_{\xi(0)}(\xi)$, 
where $\varPhi_{\xi(0)}(\xi) := \dfrac{\xi(0) - \xi}{1 - \overline{\xi(0)}\xi}$ 
is the standard M\"obius involution of the disk. Then $u_4 \in SH_{\_}(\mathbb{D})$. 
Set $D_1 := D\big(0, C_5|\varepsilon|^2\big)$ $\subset \varPhi_{\xi(0)}^{-1}(D_0)$. 
For any $\eta \in \overline{D_1}$, let $\xi = \varPhi_{\xi(0)}(\eta) \in \overline{D_0}$. Then \eqref{inequa10} implies $u_4(\eta) = u_3(\xi) \leqslant \log \bigg|\dfrac{\varepsilon}{\delta}\bigg|+C_4 $. 

By the three-circle property for subharmonic functions, for $\xi \in \overline{\mathbb{D}} \backslash D_1$, we have
\begin{equation} \label{inequa11}
\begin{aligned}
u_4(\xi) &\leqslant \bigg(\log \bigg|\dfrac{\varepsilon}{\delta}\bigg| +C_3\bigg)
\cdot  \dfrac{\log |\xi|}{\log \big(C_4 |\varepsilon|^2\big)}\\
&=\log |\xi| \cdot \bigg(\dfrac{1 - \dfrac{\log |\delta|}{\log |\varepsilon|}+\dfrac{C_3}{\log |\varepsilon|}}{2+\dfrac{\log C_4}{\log |\varepsilon|}}\bigg) .
\end{aligned}
\end{equation}
For any $\xi \in \overline{\mathbb{D}} \backslash D_0$, $\varPhi_{\xi(0)}^{-1}(\xi) = \varPhi_{\xi(0)}(\xi) \subset \overline{\mathbb{D}} \backslash D_1$. From \eqref{inequa11}, we get
\begin{multline*} 
u_2(\xi) =u_4\big(\varPhi_{\xi(0)}(\xi)\big) + \log |\xi| \\
\leqslant \log |\xi| + \log |\varPhi_{\xi(0)}(\xi)|\cdot\bigg(\dfrac{1 - \dfrac{\log |\delta|}{\log |\varepsilon|}+\dfrac{C_3}{\log |\varepsilon|}}{2+\dfrac{\log C_4}{\log |\varepsilon|}}\bigg).
\end{multline*}
On the other hand, for $\varepsilon$ small enough, $\bigg|\dfrac{\|z\|_\infty}{r_\varepsilon} - \xi(0)\bigg| \geqslant r_0 = C|\varepsilon|^2$. Pick $\xi = \dfrac{\|z\|_\infty}{r_\varepsilon} \in \mathbb{D} \backslash D_0$, then
\begin{equation*} 
\begin{aligned}
u(z_1, z_2) &= u\bigg(\Psi\bigg(\dfrac{\|z\|_\infty}{r_\varepsilon}\bigg)\bigg) = u_2\bigg(\dfrac{\|z\|_\infty}{r_\varepsilon}\bigg)\\
&\leqslant \log \bigg(\dfrac{\|z\|_\infty}{r_\varepsilon}\bigg) + \log \bigg|\varPhi_{\xi(0)}\bigg(\dfrac{\|z\|_\infty}{r_\varepsilon}\bigg)\bigg|.\bigg(\dfrac{1 - \dfrac{\log |\delta|}{\log |\varepsilon|}+\dfrac{C_3}{\log |\varepsilon|}}{2+\dfrac{\log C_4}{\log |\varepsilon|}}\bigg).
\end{aligned}
\end{equation*}
Letting $\varepsilon$ tend to $0$, 
 \begin{equation*}
\underset{\varepsilon \longrightarrow 0}{\lim} \log \bigg|\varPhi_{\xi(0)}\bigg(\dfrac{\|z\|_\infty}{r_\varepsilon}\bigg)\bigg| = \log \|z\|_\infty,
 \end{equation*}
for any $z = (z_1, z_2) \in \mathbb{D}^2 \backslash \{z_2=0\}$. 

We now use the hypothesis $\underset{\varepsilon \longrightarrow 0}{\lim} \dfrac{\log |\delta|}{\log|\varepsilon|} = 0$, and get $\underset{\varepsilon \longrightarrow 0}{\lim \sup}\; u(z) \leqslant \log \|z\|_\infty + \dfrac{1}{2}\log \|z\|_\infty = \dfrac{3}{2}\log \|z\|_\infty$, for any $z = (z_1, z_2) \in \mathbb{D}^2 \backslash \{z_2=0\}$.
\hfill $\square$

\begin{proof*}{\it Proof of Lemma \ref{smalldisk}.}

First we define a family of analytic disks going through $a_0^\eps = (0,0)$
and $a_2^\eps=(\rho,\delta\rho)$.
For $\zeta \in \mathbb{D}$, let
 \begin{equation*} 
\varphi_\lambda(\zeta) := \big(\zeta, \delta\zeta + \lambda\zeta(\zeta-\rho)\big),
\end{equation*}
where 
\begin{equation} 
\label{inequa7}
|\lambda| \leqslant \dfrac{1-|\delta|}{1+|\rho|} < 1.
\end{equation} 
We have $\varphi_\lambda(\zeta) \in \mathbb{D}^2$ for any $\zeta \in \mathbb{D}$
because of
\begin{equation*} 
\big|\delta\zeta + \lambda\zeta(\zeta-\rho)\big| \leqslant |\delta| + |\lambda|\big(1+|\rho|\big) \leqslant |\delta| + \dfrac{1-|\delta|}{1+|\rho|}.\big(1+|\rho|\big) = 1.
\end{equation*}

Let $u_1(\zeta) :=  u \circ \varphi_\lambda(\zeta) = u\big(\zeta, \delta\zeta+\lambda\zeta(\zeta-\rho)\big)$. 
From the definition of the Green function, we have
$$
u_1(\zeta) \leqslant \log\max \big(|\zeta|, |\zeta|\big|\delta+\lambda(\zeta-\rho)\big|\big)+O(1) = \log |\zeta|+O(1),
$$
because $\big|\delta+\lambda(\zeta-\rho)\big| \leqslant |\delta| + |\lambda|\big(1+|\rho|\big) \leqslant 1$.
Furthermore, 
\begin{multline*}
u_1(\zeta) = u\big(\rho+\zeta-\rho, \delta\rho+\delta(\zeta-\rho)+\lambda\zeta(\zeta-\rho)\big) \\
\leqslant \log\max \big(\big|\zeta-\rho\big|, \big|\zeta-\rho\big| \cdot \big|\delta+\lambda\zeta\big|\big)+O(1) \leqslant \log \big|\zeta-\rho\big|+O(1).
\end{multline*}
Therefore for $\zeta \in \mathbb{D}$ 
\begin{equation} \label{inequa6}
u_1(\zeta) \leqslant G^{\mathbb{D}}_{\{0,\rho\}}(\zeta) = \log \big|\zeta\dfrac{\rho-\zeta}{1-\overline{\rho}\zeta}\big|.
\end{equation} 

This provides a certain upper bound for the values of $u$ on the union of the
ranges of the disks $\varphi_\lambda$.  We want to see how it will affect $u_2$,
the restriction of $u$ on the straight disk $\Psi(\mathbb{D})$ going through
$a_1^\eps=(\eps,0)$ and $z$.

We look for $\zeta = \zeta(\lambda) \in \mathbb{D}$ such that $\varphi_\lambda(\zeta)
=\Psi(\xi)  \in \Psi(\mathbb{D})$.  Then $(\varepsilon+Z_1\xi, Z_2\xi) = \big(\zeta, \delta\zeta+\lambda\zeta(\zeta-\rho)\big)$, thus
$\zeta = \varepsilon+Z_1\xi$
and substituting into the equation for the second coordinates,
\begin{equation} \label{inequa9}
Z_2\xi = \delta(\varepsilon+Z_1\xi) + \lambda(\varepsilon+Z_1\xi)(\varepsilon - \rho + Z_1\xi).
\end{equation}
For $z_2\neq 0$, $Z_2 - \delta Z_1 \not= 0$ for $|\varepsilon| < \varepsilon_0$. 
If $\lambda = 0$ the solution of \eqref{inequa9} is
 $$
 \xi(0) := \dfrac{\delta\varepsilon}{Z_2 - \delta Z_1}. 
 $$
If $\lambda \not= 0$, let us write the solution of \eqref{inequa9}) in the following form: $\xi = \xi(\lambda) = \xi(0) + \beta(\lambda) = \dfrac{\delta\varepsilon}{Z_2 - \delta Z_1} + \beta(\lambda)$; $\beta := \beta(\lambda)$. Then \eqref{inequa9} transforms into
\begin{equation} \label{inequa090}
(Z_2-\delta Z_1)\beta = \lambda\bigg(\varepsilon+\dfrac{Z_1\delta\varepsilon}{Z_2-\delta Z_1} + Z_1 \beta\bigg)\bigg(\varepsilon - \rho + \dfrac{Z_1\delta\varepsilon}{Z_2 - \delta Z_1} + Z_1 \beta\bigg).
\end{equation}
This equation is of the form
\begin{equation}
\label{inequa09} 
a \beta^2 + (\theta_1 -b_0) \beta + c =0
\end{equation}
where  $a := \lambda Z_1^2$, $b_0 = Z_2 - \delta Z_1$, $\theta_1 = O(\varepsilon)$, $c = O(\varepsilon^2)$. The solutions of \eqref{inequa09} are
\begin{equation*} 
\begin{aligned}
\dfrac{-b \pm \sqrt{b^2 - 4ac}}{2a} &= \dfrac{\big(b_0 - \theta_1\big) \pm \sqrt{b_0^2 - 2b_0\theta_1 + \theta_1^2 - 4ac}}{2a}\\
&= \dfrac{b_0}{2a}\bigg[1 - \dfrac{\theta_1}{b_0} \pm \sqrt{1-\dfrac{2\theta_1}{b_0}+\dfrac{\theta_1^2}{b_0^2}-\dfrac{4ac}{b_0^2}}\bigg].
\end{aligned}
\end{equation*}
One of them satisfies
\begin{equation*}
\beta(\lambda) = \dfrac{b_0}{2a}\bigg[O\bigg(\dfrac{\theta_1^2}{2b_0^2}-\dfrac{2ac}{b_0^2}\bigg)\bigg] = O(\varepsilon^2).
\end{equation*}

On the other hand, \eqref{inequa090} implies 
\begin{equation*}
\lambda = \dfrac{\beta\big(Z_2-\delta Z_1\big)}{\bigg(\varepsilon+\dfrac{Z_1\delta\varepsilon}{Z_2-\delta Z_1} + Z_1 \beta\bigg)\bigg(\varepsilon - \rho + \dfrac{Z_1\delta\varepsilon}{Z_2 - \delta Z_1} + Z_1 \beta\bigg)}.
\end{equation*}
Since $|\varepsilon - \rho| \geqslant \dfrac{1}{2}|\varepsilon|$ and $\big|Z_2-\delta Z_1\big| \leqslant 1+|\delta|$, it follows that
\begin{equation*}
|\lambda| = \dfrac{\big|(Z_2-\delta Z_1)\big||\beta|}{\bigg|\big(\varepsilon + O(|\delta\varepsilon|)\big)\big(\varepsilon-\rho + O(|\delta\varepsilon|)\big)\bigg|} \leqslant \dfrac{|\beta|\big(1+|\delta|\big)}{1/2|\varepsilon|^2}.
\end{equation*}
So $\dfrac{|\beta|\big(1+|\delta|\big)}{1/2|\varepsilon|^2} \leqslant \dfrac{1-|\delta|}{1+|\rho|}$, i.e.
\begin{equation*} 
|\beta(\lambda)| \leqslant \dfrac{1-|\delta|}{1+|\rho|} \cdot \dfrac{1}{1+|\delta|}
\cdot \dfrac{1}{2}|\varepsilon|^2 < \dfrac{1}{2}|\varepsilon|^2,
\end{equation*}
and $|\lambda| \leqslant \dfrac{1-|\delta|}{1+|\rho|}$.
 
Set $z_0 := \xi(0) = \dfrac{\delta\varepsilon}{Z_2-\delta Z_1}$.
For $|\eps|$, therefore $|\delta|$ small enough, $|z_0| = \bigg|\dfrac{\delta\varepsilon}{Z_2-\delta Z_1}\bigg| \leqslant C_0(z)|\varepsilon| < 1$. Now consider a disc $D_0 := D(z_0;r_0)$, where $r_0 := \frac12 |\varepsilon|^2$. 

Since $\xi \in \overline{D_0}$, $|\zeta(\xi)| = |\varepsilon + Z_1\xi| \leqslant |\varepsilon| + |\xi| \leqslant (|\varepsilon| + C_0|\varepsilon| + \frac12 |\varepsilon|^2) \leqslant C_1|\varepsilon|$, where $C_0 := C_0(z)$ and $0 < C_1$ do not depend on $\varepsilon$. So
\begin{equation}\label{inequa091} 
\begin{aligned}
\bigg|(\varepsilon+Z_1 \xi)\dfrac{\rho-\varepsilon-Z_1\xi}{1-\overline{\rho}(\varepsilon+Z_1 \xi)}\bigg| &\leqslant C_1|\varepsilon| .\dfrac{1/2|\varepsilon| + C_1|\varepsilon|}{\big|1-1/2C_1|\varepsilon|^2\big|}\\
&\leqslant C_1\dfrac{1/2+C_1}{|1-1/2C_1|}|\varepsilon|^2 \leqslant C_2|\varepsilon|^2,
\end{aligned}
\end{equation}
where $0 < C_2$ does not depend on $\varepsilon$. So for $|\varepsilon|$ small enough and at least $ \leqslant \dfrac{1}{C_1}$,  there is $\zeta(\xi) \in \mathbb{D}$ 
such that $\Psi(\xi) = \varphi_\lambda\big(\zeta(\xi)\big)$. 
Thus
$$
u_2(\xi) = u\big(\Psi(\xi)\big) = u\big(\varphi_\lambda\big(\zeta(\xi)\big)\big) = u_1(\zeta(\xi)) = u_1(\varepsilon+Z_1 \xi). 
$$
From (\ref{inequa6}) and (\ref{inequa091}), we deduce
\begin{equation*} 
\begin{aligned}
u_1(\varepsilon+Z_1 \xi) &\leqslant G^{\mathbb{D}}_{\{0, \rho\}}(\varepsilon+Z_1 \xi) = \log \bigg|(\varepsilon+Z_1 \xi)\dfrac{\rho-\varepsilon-Z_1\xi}{1-\overline{\rho}(\varepsilon+Z_1 \xi)}\bigg| \\
&\leqslant \log |\varepsilon|^2 + \log C_2.
\end{aligned}
\end{equation*}
For any $\xi \in \overline{D_0}$, 
\begin{equation*} 
||z_0||  = \dfrac{\big|\delta\varepsilon\big|}{\big|Z_2-\delta Z_1\big|}> \dfrac{|\delta||\varepsilon|}{\big|Z_2-\delta Z_1\big|} > \dfrac{|\delta||\varepsilon|}{1+|\delta|} > \dfrac{1}{2}|\varepsilon|^2 = r_0,
\end{equation*} 
for $|\eps|$ small enough depending on $z$, so $\xi \not= 0$ and 
$$
u_3(\xi) = u_2(\xi) - \log |\xi| \leqslant \log |\varepsilon|^2 - \log |\xi| + \log C_2.
$$
In addition, $\big|\xi-\xi(0)\big| = |\beta(\lambda)| < 1/2|\varepsilon|^2$ and $|\varepsilon\delta| \gg |\varepsilon|^2$ imply that 
$$|\xi| \geqslant \dfrac{1}{\big|Z_2-\delta Z_1\big|}|\varepsilon\delta| - \dfrac{1}{2} |\varepsilon|^2 > \bigg(\dfrac{1}{\big|Z_2-\delta Z_1\big|} - \dfrac{1}{2}\bigg)|\delta\varepsilon| >0
$$ 
(because $\big|Z_2-\delta Z_1\big| < 1 + |\delta| < 2$). Then 
$$
u_3(\xi) \leqslant \log |\varepsilon|^2 - \log |\varepsilon\delta| - \log \bigg(\dfrac{1}{\big|Z_2-\delta Z_1\big|} - \dfrac{1}{2}\bigg) + \log C_2,
$$
for any $\xi \in \overline{D_0}$. Finally 
\begin{equation} \label{inequa10}
u_3(\xi) \leqslant \log \bigg|\dfrac{\varepsilon}{\delta}\bigg| + \; C_3,  \mbox{ with } C_3 := C_3(z) = - \log \bigg(\dfrac{1}{\big|Z_2-\delta Z_1\big|} - \dfrac{1}{2}\bigg) + \; \log C_2.
\end{equation}

\end{proof*}

{}

\end{document}